\documentclass[12pt]{amsart}
\usepackage{amsthm, amssymb}
\usepackage{amsfonts, amscd}
\usepackage{epsfig,multicol}
\theoremstyle{plain} \numberwithin{equation}{section}
\newtheorem{theo}{Theorem}[section]

\newtheorem{lemm}[theo]{Lemma}

\theoremstyle{definition}

\newtheorem*{rema}{Remark}
\newtheorem*{theoap}{Theorem}
\newtheorem*{propap}{Proposition}

\def\Z{\mathbb Z}
\def\C{\mathbb C}
\def\R{\mathbb R}

\def\s{\sigma}
\def\GG{\mathcal G}
\def\G{\Gamma}
\def\e{\mathcal E}
\def\H{\Gamma}

\def\f{\rho}
\def\g{\rho}
\def\n{\ell}
\def\Phi{\bar \rho}
\def\A{\phi_A}
\def\B{\phi_B}
\def\Q{Q}
\def\tQ{\tilde Q}
\def\r{\bar r}
\def\TR{T^n(\R)}
\def\XR{X(\R)}

\newcommand{\ZZ}{{\mathbb Z}}
\newcommand{\al}{\alpha}
\newcommand{\be}{\beta}
\newcommand{\ra}{{\rightarrow}}

\DeclareMathOperator{\GL}{GL}
\DeclareMathOperator{\Aut}{Aut}

\DeclareMathOperator{\Ker}{Ker}

\begin{document}
\title{Cohomological rigidity of real Bott manifolds}
\author[Y. Kamishima and M. Masuda]{Yoshinobu Kamishima and Mikiya Masuda}
\address{Department of Mathematics, Tokyo Metropolitan University, 
Minami-Ohsawa 1-1, Hachioji, Tokyo 192-0397, Japan} 
\email{kami@tmu.ac.jp}

\address{Department of Mathematics, Osaka City
University, Sumiyoshi-ku, Osaka 558-8585, Japan.}
\email{masuda@sci.osaka-cu.ac.jp}

\date{\today}
\thanks{The first author was partially supported by Grant-in-Aid for
Scientific Research 2034001 
and the second author was partially supported by Grant-in-Aid for
Scientific Research 19204007}
%\subjclass[2000]{Primary 57S15, 14M25; Secondary 57S25}

\begin{abstract}
A real Bott manifold is the total space of iterated $\R P^1$ bundles starting with 
a point, where each $\R P^1$ bundle is projectivization of a Whitney sum of two real line bundles.  
We prove that two real Bott manifolds are diffeomorphic if their cohomology rings 
with $\Z/2$ coefficients are isomorphic.  

A real Bott manifold is a real toric manifold and admits a flat riemannian metric 
invariant under the natural action of an elementary abelian 2-group.  
We also prove that the converse is true, namely a real toric manifold  
which admits a flat 
riemannian metric invariant under the action of an elementary abelian 2-group 
is a real Bott manifold.  
\end{abstract}

\maketitle 

%\vskip .3cm

\section{Introduction}

A fundamental result in the theory of toric varieties says that 
the categories of toric varieties (over the complex numbers $\C$) and 
fans are equivalent (see \cite{oda88}).  
This reduces the classification of toric varieties to that of fans. 
Among toric varieties, compact smooth toric varieties which we call 
toric manifolds are well studied and the classification as varieties is 
completed for some classes of toric manifolds 
(see \cite{klei88}, \cite{oda88}, \cite{sato06} for example). 

However, not much is known for the topological classification of toric 
manifolds, and the following problem is addressed in \cite{ma-su07}
(see also \cite{ch-ma-su08-2}, \cite{ma-pa08}). 

\medskip
\noindent
{\bf Cohomological rigidity problem for toric manifolds.}
Are two toric manifolds diffeomorphic (or homeomorphic) 
if their cohomology 
rings with integer coefficients are isomorphic as graded rings? 

\medskip
As is well-known, there are many closed smooth manifolds which are not 
homeomorphic but have isomorphic cohomology rings.  So the problem above 
seems unlikely but no counterexample is known and there are 
some partial affirmative solutions to the problem (see \cite{ch-ma-su08-2}, 
\cite{ma-pa08}, \cite{ma-su07}). 

The set $\XR$ of real points in a toric manifold $X$ is called 
a real toric manifold. 
It appears as the fixed point set of the complex conjugation on $X$. 
For example, when $X$ is a complex projective space $\C P^n$, $\XR$ is 
a real projective space $\R P^n$.  It is known that 
\begin{equation*} \label{Z2}
H^*(\XR;\Z/2)\cong H^{2*}(X;\Z)\otimes \Z/2
\end{equation*}
for any toric manifold $X$ where $\Z$ denotes the integers and 
$\Z/2=\{0,1\}$, 
and one may ask the same question as the rigidity problem above for 
real toric manifolds with $\Z/2$ coefficients, namely 

\medskip
\noindent
{\bf Cohomological rigidity problem for real toric manifolds.}
Are two real toric manifolds diffeomorphic (or homeomorphic) if their 
cohomology rings with $\Z/2$ coefficients are isomorphic as graded rings? 

\medskip
In this paper we are concerned with a sequence of $\R P^1$ bundles 
\begin{equation} \label{tower}
M_n\stackrel{\R P^1}\longrightarrow M_{n-1}\stackrel{\R P^1}\longrightarrow 
\cdots\stackrel{\R P^1}\longrightarrow M_1
\stackrel{\R P^1}\longrightarrow M_0=\{\textrm{a point}\}
\end{equation}
such that $M_i\to M_{i-1}$ for $i=1,\dots,n$ is the 
projective bundle of a Whitney sum of 
two real line bundles over $M_{i-1}$, where one of the two 
line bundles may be assumed to be trivial without loss of generality. 
Grossberg-Karshon \cite{gr-ka94} considered the sequence above in the 
complex case 
and named it a {\em Bott tower} of height $n$.  Following them, we call the 
sequence above a {\em real Bott tower} of height $n$.  
The top manifold $M_n$ of a real Bott tower is a real toric manifold. 
We call it a {\em real Bott manifold}. 
The main purpose of this paper is to prove the following which answers 
the cohomological rigidity problem affirmatively for real Bott manifolds. 

\begin{theo} \label{main}
Two real Bott manifolds are diffeomorphic if their cohomology 
rings with $\Z/2$ coefficients are isomorphic as graded rings. 
\end{theo}

Although real toric manifolds have similar properties to 
toric manifolds, there is one major difference, that is, 
a real toric 
manifold is not simply connected while a toric manifold is 
simply connected. 
Real toric manifolds provide many examples of aspherical 
manifolds and real Bott manifolds are examples of flat riemannian 
manifolds. In fact, any real toric manifold of dimension $n$ 
supports an action of an elementary abelian $2$-group $\TR$ of rank $n$ 
and real Bott manifolds of dimension $n$ admit a flat riemannian metric 
invariant under the action of $\TR$. The following shows that these are the 
only examples among real toric manifolds. 

\begin{theo} \label{main2}
A real toric manifold of dimension $n$ which 
admits a flat riemannian metric invariant under the action of $\TR$ 
is a real Bott manifold. 
\end{theo}

This paper is organized as follows.  We describe the cohomology ring 
and the fundamental group of a real Bott manifold in Sections~\ref{sect2} 
and ~\ref{fund}.  In Section~\ref{coho} we find necessary 
and sufficient conditions for an isomorphism between 
cohomology rings of real Bott manifolds to satisfy in terms of matrices.  
Using the conditions, 
we construct a monomorphism between the fundamental groups of 
the real Bott manifolds in Section~\ref{sect5}.  It 
may not be an isomorphism but the existence of the monomorphism implies that 
the fundamental groups are isomorphic, which 
is done in Section~\ref{sect6} by studying group extensions. 
Since real Bott manifolds are flat riemannian manifolds, the isomorphism 
of the fundamental groups implies Theorem~\ref{main} by a theorem 
of Bieberbach.  
In Section~\ref{sect7} we enumerate diffeomorphism classes 
in real Bott manifolds of dimension up to $4$. This result is obtained 
in \cite{nazr08} independently by a different method. 
Theorem~\ref{main2} is proved in Section~\ref{sect8}. 
In Section~\ref{sect9} we view real Bott manifolds 
from a viewpoint of small covers introduced in \cite{da-ja91}. 
In the Appendix, we give a proof on a (probably known) fact used in 
Section~\ref{sect6}.

\section{Cohomology rings} \label{sect2}

We shall describe the cohomology ring
of the real Bott manifold $M_n$ in the tower (\ref{tower}). 

We recall a general well-known fact. 
Let $E\to X$ be a real vector bundle of rank $m$ over a topological space 
$X$ and let $P(E)$ be the projectivization of $E$.  As is well-known, 
$H^*(P(E);\Z/2)$ is an algebra over $H^*(X)$ through the projection 
map from $P(E)$ to $X$ and the algebra structure is described as 
\begin{equation} \label{proj}
H^*(P(E);\Z/2)=H^*(X;\Z/2)[x]\Big/(\sum_{i=0}^m w_i(E)x^{m-i})
\end{equation}
where $w_i(E)$ denotes the $i$-th Stiefel-Whitney class of $E$ and 
$x$ is given by the first Stiefel-Whitney class of the canonical line 
bundle over $P(E)$. 
Moreover, the Stiefel-Whitney class of $T_fP(E)$ the tangent bundle 
along the fibers of $P(E)$ is given by 
\[
w(T_f(P(E)))=\sum_{i=0}^mw_i(E)(1+x)^{m-i},
\]
in particular, 
\begin{equation} \label{w1Tf}
w_1(T_f(P(E)))=w_1(E) 
\end{equation}
when $m$ is even.

Now we return to the tower (\ref{tower}).  
By definition $M_j=P(L_{j-1}\oplus\underline{\R})$ 
with some line bundle $L_{j-1}$ over $M_{j-1}$ for $j=1,\dots,n$ , 
where $\underline{\R}$ 
denotes the trivial line bundle.  Let $\gamma_j$ be the canonical line 
bundle over $M_j$ and set $x_j=w_1(\gamma_j)$.  We use the same notation 
$\gamma_j$ (resp. $x_j$) for the pullback of $\gamma_j$ (resp. $x_j$) 
by compositions of 
projections $M_k\to M_{k-1}\to \dots \to M_j$ where $k>j$.  
Then the repeated use of (\ref{proj}) shows 
\begin{equation} \label{bcoho}
H^*(M_k;\Z/2)=\Z/2[x_1,\dots,x_k]\big/\big(x_j(x_j+w_1(L_{j-1}))\mid 
j=1,\dots,k\big).
\end{equation}

Since $H^1(M_{j-1};\Z/2)$ is additively generated by $x_1,\dots,x_{j-1}$ 
and $L_{j-1}$ is a line bundle over $M_{j-1}$, one can uniquely write 
\begin{equation} \label{w1}
w_1(L_{j-1})=\sum_{i=1}^{j-1}A^i_jx_i\quad\text{with $A^i_j\in \Z/2$}
\end{equation}
where $j=2,\dots,n$.  
As is well-known, line bundles are classified by 
their first Stiefel-Whitney classes and 
the first Stiefel-Whitney class behaves 
additively for tensor products of line bundles; so it follows from 
\eqref{w1} that 
\begin{equation} \label{L}
L_{j-1}=\gamma_1^{A^1_j}\otimes \dots\otimes \gamma_{j-1}^{A^{j-1}_j}.
\end{equation}
For convenience, we set $A^i_j=0$ unless $i<j$ and form 
a square matrix $A$ of size $n$ with $A^i_j$ as an $(i,j)$ entry. 
$A$ is an upper triangular matrix with zero diagonal entries. 

The observation above implies that the tower (\ref{tower}) is completely 
determined by the matrix $A$.  So we may denote $M_n$ by $M(A)$. 
For later use we record the ring structure of $H^*(M(A);\Z/2)$ 
as a lemma which follows from (\ref{bcoho}) and (\ref{w1}). 

\begin{lemm} \label{HBn}
Let $A$ and $M(A)$ be as above.  
Then $H^*(M(A);\Z/2)$ is generated by degree one 
elements $x_1,\dots,x_n$ as a graded ring with $n$ relations 
\[
x_j^2=x_j\sum_{i=1}^n A^i_jx_i \quad\text{for $j=1,\dots,n$}.
\]
\end{lemm}

We conclude this section with the following lemma. 

\begin{lemm} \label{orien}
The real Bott manifold $M(A)$ is orientable if and only if the sum of 
entries is zero in $\Z/2$ for each row of $A$.
\end{lemm}

\begin{proof}
The repeated use of \eqref{w1Tf} together with \eqref{w1} shows that 
\[
\begin{split}
w_1(M(A))&=\sum_{j=1}^n w_1(L_{j-1}\oplus\underline{\R})
=\sum_{j=1}^n w_1(L_{j-1})\\
&=\sum_{j=1}^n \sum_{i=1}^{j-1}A^i_jx_i
=\sum_{i=1}^n (\sum_{j=1}^n A^i_j)x_i.
\end{split}
\]
Since $M(A)$ is orientable if and only if $w_1(M(A))=0$, the lemma follows 
from the identity above. 
\end{proof}

\section{Fundamental groups} \label{fund}

A general description of the fundamental group of 
an arbitrary real toric manifold is given in \cite{uma04} motivated by 
the work \cite{da-ja91}. 
In this section, we shall describe the fundamental group of $M(A)$ in 
a direct way. 

Let $s_i$ $(i=1,\dots,n)$ 
be an Euclidean motion on $\R^n$ defined by 
\begin{equation} \label{si}
\begin{split}
s_i(u_1,\dots,u_n)&=(u_1,\dots,u_{i-1}, u_i+\frac{1}{2}, 
(-1)^{A_{i+1}^i}u_{i+1},\dots, (-1)^{A_{n}^i}u_n)\\
&=( (-1)^{A_1^i}u_1,\dots, (-1)^{A_n^i}u_n)+\frac{1}{2}e^i
\end{split}
\end{equation}
where $e^1,\dots,e^n$ denote the standard basis of $\R^n$. 
The group $\G(A)$ generated by $s_1,\dots,s_n$ is a crystallographic 
group. In fact, the subgroup generated by $s_1^2,\dots,s_n^2$ 
consists of all translations by $\Z^n$.  The action of $\G(A)$ on $\R^n$ 
is free and the orbit space $\R^n/\G(A)$ is compact.  

\begin{lemm} \label{rflat}
$\R^n/\G(A)$ is diffeomorphic to $M(A)$. Therefore $M(A)$ is a riemannian flat 
manifold with $\G(A)$ as the fundamental group. 
\end{lemm}

\begin{proof}
Let $\G_k$ $(k=1,\dots,n)$ 
be a subgroup of $\G(A)$ generated by $s_1,\dots,s_k$. 
It acts on $\R^k$ by restricting the action of $\G(A)$ on $\R^n$.  
We claim that a sequence of projections 
$$\R^n/\G_n\to \R^{n-1}/\G_{n-1}\to\dots \to\R^1/\G_1\to \{0\}$$
agrees with the real Bott tower (\ref{tower}).  
The lemma follows from the claim. 

We shall prove the claim by induction on height.  
It is obviously true up to height one.  Suppose it is true up to height 
$j-1$. We note that the line bundle $\gamma_{i}$ over $M_{j-1}$ for 
$i\le j-1$ is obtained as 
the quotient of $\R^{j-1}\times\R$ by the diagonal action of $\G_{j-1}$ 
where the action of $\G_{j-1}$ on the second factor 
$\R$ is given through a homomorphism 
$\G_{j-1}\to \{\pm 1\}$ sending $s_i$ to $-1$ and the others $s_\ell$ 
($\ell\not=i$) to $1$.  
This together with (\ref{L}) shows that the line bundle $L_{j-1}$ in 
(\ref{L}) 
is obtained as the quotient of $\R^{j-1}\times \R$ by the diagonal action 
of $\G_{j-1}$ where the action of $\G_{j-1}$ on the second factor 
$\R$ is given through a homomorphism $\G_{j-1}\to \{\pm 1\}$ sending $s_i$ to 
$(-1)^{A^i_{j}}$ for $i\le j-1$.  Therefore the action of $\G_{j-1}$ on 
$\R^{j-1}\times \R=\R^{j}$ 
is nothing but the restriction of the action of $\G_{j}$ to $\G_{j-1}$ 
while the action of $s_{j}$ on $\R^{j}$ is trivial on the first 
$(j-1)$ coordinates and translation by $1/2$ on the last coordinate. 

We consider a map  
\begin{equation*}
\begin{split}
\R^j=\R^{j-1}\times\R\ &\ra\  
(\R^{j-1}\times\R)/\G_{j-1}\oplus\underline{\R}=
L_{j-1}\oplus\underline{\R}; \\
(x,\theta)\ &\mapsto\ \big([x,\sin2\pi\theta],\cos2\pi\theta\big).
\end{split}
\end{equation*}
Since $s_i(x,\theta)=(s_ix,(-1)^{A^i_j}\theta)$ for $i\le j-1$ and $s_j^2(x,\theta)=(x,\theta+1)$, 
the map above 
is invariant under the action of $\G_{j-1}$ and $s_j^2$ and factors through  
a diffeomorphism from the orbit space $\R^j/\langle \G_{j-1},s_j^2\rangle$ 
onto the unit circle bundle of $L_{j-1}\oplus\underline{\R}$. 
Furthermore, since $\G_j=\langle \G_{j-1},s_j\rangle$ and $s_j(x,\theta)=(x,\theta+\frac{1}{2})$, 
the map induces a diffeomorphism from $\R^j/\G_j$ onto 
the projectivization $P(L_{j-1}\oplus\underline{\R})=M_j$. 
This shows that the projection $\R^{j}/\G_{j}\to \R^{j-1}/\G_{j-1}$ 
agrees with the projection $M_j\to M_{j-1}$, completing the induction step. 
\end{proof}

We shall investigate the structure of $\G(A)$. 

\begin{lemm} \label{lemm1}
For $i<\ell$, $s_\ell s_i=s_is_\ell^{(-1)^{A_\ell^i}}$, i.e. 
\[
s_\ell s_i=\begin{cases} s_is_\ell^{-1} \quad&\text{if $A_\ell^i=1$,}\\
s_is_\ell\quad&\text{if $A_\ell^i=0$.}
\end{cases}
\]
\end{lemm}

\begin{proof} Easy to check. 
%It is easy to see that the images of $x\in \R^n$ 
%by $s_is_j$ and $s_js_i^{(-1)^{{B}_i^j}}$ agree except the $j$-th component 
%regardless of the value of $B_i^j$.   
\end{proof}

\begin{lemm} \label{lemm2}
Let $\GG(A)$ be the group generated by $\s_1,\dots,\s_n$ with the relations 
in Lemma~\ref{lemm1} for $\s_j$'s instead of $s_j$'s. Then the homomorphism 
$\psi\colon \GG(A)\to \G(A)$ defined by $\psi(\s_j)=s_j$ for $j=1,\dots,n$ 
is an isomorphism.
\end{lemm}

\begin{proof} 
Using the relations, one can express an element $\s$ of $\GG(A)$ as 
$\s_1^{a_1}\s_2^{a_2}\dots \s_n^{a_n}$ with $a_1,\dots,a_n\in\Z$. 
Suppose $\psi(\s)=s_1^{a_1}s_2^{a_2}\dots s_n^{a_n}$ is 
the identity element.  Then $\psi(\s)$ fixes 
any element of $\R^n$.  But it maps the origin of $\R^n$ to 
$\frac{1}{2}\sum_{j=1}^n\epsilon_ja_je_j$, where $\epsilon_j=\pm 1$, 
and the image must again be the origin, so 
we have $a_j=0$ for any $j$.  This shows that $\s$ is the identity and 
$\psi$ is injective.  The surjectivity of $\psi$ is trivial. 
\end{proof}

$s_j^2$ is a translation of $\R^n$ by $e_j$ so that $s_1^2,\dots,s_n^2$ 
commute with each other and generate a free abelian subgroup $N$. 
The images of $s_j$'s in the quotient $\G(A)/N$ commute with each other, 
which easily follows from Lemma~\ref{lemm1}, so that $\G(A)/N$ is an 
elementary abelian $2$-group.  
We identify $N$ with $\Z^n$ and $\G(A)/N$ with $(\Z_2)^n$ in a natural way 
where $\Z_2=\{\pm 1\}$ 
and obtain a short exact sequence:
\begin{equation} \label{exact}
0\to \Z^n \to \G(A)\to (\Z_2)^n\to 1.
\end{equation}

One may think of $M(A)=\R^n/\G(A)$ as the orbit space of the torus 
$\R^n/\Z^n$ by the induced action of $\G(A)/\Z^n=(\Z_2)^n$. 
We shall explicitly describe the action using complex numbers $\C$. 
Let $S^1$ denote the unit circle of $\C$. 
We identify $\R/\Z$ with $S^1$ (and hence $\R^n/\Z^n$ with $(S^1)^n$) 
through the exponential map sending 
$u\in \R$ to $\exp(2\pi \sqrt{-1} u)\in \C$. 
For $z\in S^1$ and $a\in \Z/2=\{0,1\}$ we define 
\[
z(a):=\begin{cases} z \quad&\text{if $a=0$,}\\
\bar z\quad&\text{if $a=1$.}
\end{cases}
\]
Then the induced action of $s_i$ defined in \eqref{si} on $(S^1)^n$ is 
given by 
\[
(z_1,\dots,z_n)\to (z_1,\dots,z_{i-1},-z_i,z_{i+1}(A^i_{i+1}),\dots,
z_n(A^i_n)). 
\]

\section{An isomorphism between cohomology rings} \label{coho}

As is described in Lemma~\ref{HBn}, 
$H^*(M(A);\Z/2)=R_1$ is a graded algebra over $\Z/2$ 
generated by degree one elements $x_1,\dots,x_n$ with $n$ relations
\begin{equation} \label{rel}
x_j^2=x_j\sum_{i=1}^n A_j^ix_i \quad(j=1,\dots,n).
\end{equation}
The set $V_1$ of degree one elements in $R_1$ with vanishing squares 
forms a vector space over $\Z/2$ of positive dimension.  Set $n_1=\dim V_1$. 
Permuting the suffixes of $x_1,\dots,x_n$, we may assume that 
the first $n_1$ elements $x_1,\dots,x_{n_1}$ form a basis of $V_1$. 
We consider the quotient graded ring $R_2=R_1/(V_1)$ where $(V_1)$ denotes 
the ideal in $R_1$ generated by $V_1$. 
Similarly, the set $V_2$ of degree one elements in $R_2$ with vanishing 
squares forms a vector space over $\Z/2$ of positive dimension. 
Set $n_2=\dim V_2$.  Permuting the suffixes 
of $x_{n_1+1},\dots,x_n$, we may assume that the image of 
$x_{n_1+1},\dots,x_{n_1+n_2}$ in the quotient ring $R_2$ forms a basis of 
$V_2$. Then consider the quotient graded ring $R_3=R_2/(V_2)$ and repeat 
the same argument, and so on. 
This procedure will terminate at a finite steps, say $q$ steps, so 
that we obtain a sequence of natural numbers $(n_1,\dots,n_q)$, which is an 
invariant of the cohomology ring.  We call this sequence the \emph{type} of 
$A$ or $H^*(M(A);\Z/2)$. 
The argument above shows that through a suitable permutation of suffixes of 
$x_1,\dots,x_n$ we may assume that the upper triangular matrix $A$ 
decomposed into $q\times q$ blocks according to the type $(n_1,\dots,n_q)$ 
has zero matrices of sizes $n_1,\dots,n_q$ as the diagonal $q$ blocks, 
i.e. 
\begin{equation} \label{block}
A=\left(
\begin{array}{cccc} O_{n_1} & & &\ast \\ & O_{n_2} & & \\ & & \ddots& \\ 
0&  && O_{n_q} 
\end{array}\right)
\end{equation}
where $O_m$ denotes the zero matrix of size $m$ 
and $(i,i+1)$-block is non-zero for each $i=1,\dots,q-1$.
We note that permuting suffixes of $x_1,\dots,x_n$ corresponds to conjugating 
the matrix $A$ by a permutation matrix. 

Let $B$ be an upper triangular matrix of the same type and same form as 
(\ref{block}) and 
let $$\varphi\colon H^*(M(A);\Z/2)\to H^*(M(B);\Z/2)$$ be an isomorphism 
as graded rings.  We denote by $y_1,\dots,y_n$ the 
generators of $H^*(M(B);\Z/2)$.  
Since $\varphi({x_i})^2=\varphi({x_i}^2)=0$ for $1\le i\le 
n_1$, $\varphi({x_i})$ is a linear combination of $y_1,\dots,y_{n_1}$.  
In general, 
one easily sees that $\varphi(x_i)$ for $n_{j-1}+1\le i\le n_{j}$ is 
a linear combination of $y_1,\dots,y_{n_{j}}$.  This means that 
if we view $P\in \GL(n;\Z/2)$ defined by 
\begin{equation} \label{P}
(\varphi(x_1), \dots,\varphi(x_n))=(y_1,\dots,y_n)P
\end{equation}
as a $q\times q$ block matrix of type $(n_1,\dots,n_q)$, then 
$P$ is an upper triangular block matrix.  Since $P$ is non-singular, 
every diagonal block of $P$ is also non-singular.  Therefore, we may assume 
that the diagonal entries of $P$ are all one if 
necessary by permuting the suffixes of the generators $y_i$'s in each block. 
With this understood, we have 

\begin{lemm} \label{bBC}
$B=PA$ and
\[
P_j^\ell B^i_\ell=P_j^iB_j^\ell+P_j^\ell B_j^i+
P_j^\ell B_j^\ell B_\ell^i \quad\text{for $i<\ell$.}
\]
\end{lemm}

\begin{proof}
It follows from (\ref{P}) that 
\begin{equation} \label{phi}
\varphi(x_k)=\sum_{i=1}^nP_k^i y_i \quad\text{for}\quad k=1,\dots,n.
\end{equation}
We plug this in (\ref{rel}) mapped by $\varphi$ to obtain 
\begin{equation} \label{CBP}
\begin{split}
\big(\sum_{i=1}^n P_j^i y_i\big)^2&=
\big(\sum_{i=1}^n P_j^i y_i\big)\big(\sum_{k=1}^n 
\sum_{i=1}^nA_j^kP_k^iy_i\big)\\
&=\big(\sum_{i=1}^n P_j^i y_i\big)\big(\sum_{i=1}^n(PA)_j^iy_i\big)
\end{split}
\end{equation}
Comparing the coefficients of $y_iy_j$ for $i<j$ at both sides 
above and noting that $P_j^j=1$ and $(PA)_j^j=0$, we obtain
\[
B^i_j=(PA)^i_j \quad\text{for $i<j$}.
\]
(Note that the term $y_iy_j$ may appear in ${y_j}^2$ but not in 
${y_i}^2$ because $B$ is assumed to be upper triangular.) 
The identity above holds even for $i\ge j$ because the both sides then vanish. 
Therefore $B=PA$. 

More generally, comparing the coefficients of $y_iy_\ell$ for $i<\ell$ 
at the both sides of (\ref{CBP}) and replacing $PA$ by $B$, 
we obtain the latter identity in the lemma. 
\end{proof}

\section{A monomorphism between fundamental groups} \label{sect5}

Let $A, B$ and $P$ be as in Section~\ref{coho}. 
In this section we 
construct a monomorphism between the fundamental groups $\G(B)$ and $\G(A)$ 
using $P$.

Any element $s\in \G(A)$ can be expressed uniquely as 
$s=s_1^{a_1}s_2^{a_2}\dots s_n^{a_n}$ with integers $a_i$'s 
by Lemma~\ref{lemm1}.  
We denote the exponent $a_j$ of $s_j$ by $\e_j(s)$. 

\begin{lemm} \label{expo}
If $p_i,q_i\in \Z$ for $i=1,\dots,n$, then 
\[
\begin{split}
&\e_j((s_1^{p_1}s_2^{p_2}\dots s_n^{p_n})(s_1^{q_1}s_2^{q_2}\dots s_n^{q_n}))
=(-1)^{\sum_{k=1}^{j-1}q_kB_j^k}p_j+q_j\\
&\e_j((s_1^{p_1}s_2^{p_2}\dots s_n^{p_n})
(s_n^{-q_n}\dots s_2^{-q_2}s_1^{-q_1}))
=(-1)^{\sum_{k=1}^{j-1}q_kB_j^k}(p_j-q_j)
\end{split}
\]
\end{lemm}

\begin{proof}
Using Lemma~\ref{lemm1}, we see 
\begin{equation} \label{tP}
s_\ell^p s_k^q=s_k^qs_\ell^{p(-1)^{qB_\ell^k}}\quad\text{for $\ell>k$, and 
$p, q\in\Z$}
\end{equation}
and the repeated use of this identity implies the lemma. 
\end{proof}

We use notation $t_i$'s for $\G(B)$ in place of $s_i$'s for $\G(A)$.  
We regard $P$ as an \emph{integer} matrix and define 
\begin{equation} \label{f}
\f(t_r)=s_1^{P_1^r}s_2^{P_2^r}\dots s_n^{P_n^r} \quad(r=1,\dots,n).
\end{equation}
We shall check that $\f$ preserves the relations in Lemma~\ref{lemm1} 
for $\G(B)$ so that $\f$ induces a homomorphism from $\G(B)$ to $\G(A)$ by 
Lemma~\ref{lemm2}.  It follows from Lemma~\ref{expo} that 
\begin{equation} \label{1}
\e_j(\f(t_\ell t_i))=(-1)^{\sum_{k=1}^{j-1}P_k^i A_j^k}P_j^\ell+P_j^i
=(-1)^{B_j^i}P_j^\ell+P_j^i\in\Z
\end{equation}
where we used the fact $PA=B$ and $A_j^k=0$ for $k\ge j$ 
in the latter identity.  Similarly we have 
\begin{equation} \label{2}
\e_j(\f(t_i t_\ell))=(-1)^{B_j^\ell}P_j^i+P_j^\ell\in\Z
\end{equation}
and 
\begin{equation} \label{3}
\e_j(\f(t_i t_\ell^{-1}))=(-1)^{B_j^\ell}(P_j^i-P_j^\ell)\in\Z.
\end{equation}

Now suppose $i<\ell$.  When $B_\ell^i=0$, we have $t_\ell t_i=t_it_\ell$ by 
Lemma~\ref{lemm1} for $\G(B)$ and 
\[
P_j^\ell B_j^i=P_j^iB_j^\ell \in\Z/2
\] 
by Lemma~\ref{bBC}. 
An elementary case-by-case check (according to the values of $B_j^i$ 
and $B_j^\ell$) shows that the identity above ensures that 
the right hand sides at (\ref{1}) and (\ref{2}) coincide.  
When $B_\ell^i=1$, we have 
$t_\ell t_i=t_it_\ell^{-1}$ by Lemma~\ref{lemm1} for $\G(B)$ and 
\[
P_j^\ell =P_j^i B_j^\ell+P_j^\ell B_j^i+
P_j^\ell B_j^\ell\in\Z/2 \quad\text{for $i<\ell$}
\]
by Lemma~\ref{bBC}. A similar elementary case-by-case check shows that the
identity above ensures that 
the right hand sides at (\ref{1}) and (\ref{3}) coincide.  
In any case the map $\f$ preserves the relations for $\G(B)$ and $\G(A)$ and 
hence induces a homomorphism from $\G(B)$ to $\G(A)$.  

\begin{lemm} \label{odd}
The homomorphism $\f\colon \G(B)\to \G(A)$ is injective and 
\begin{enumerate}
\item[(1)] $\f(\Z^n)\subset \Z^n$ and $\Z^n/\f(\Z^n)$ is of order $\det P$ 
(which is odd),
\item[(2)] $\f$ induces an isomorphism from $\G(B)/\Z^n$ onto $\G(A)/\Z^n$. 
\end{enumerate}
Therefore 
$\f$ is an isomorphism if and only if $\det P=\pm 1$. 
\end{lemm}

\begin{proof}
It follows from (\ref{1}) with $\ell=i$ that 
\[
\e_j(\f(t_i^2))=\begin{cases} 2P_j^i \quad&\text{if $B_j^i=0$,}\\
0\quad&\text{if $B_j^i=1$.}
\end{cases}
\]
Therefore $\f$ maps the normal subgroup $\Z^n$ of $\G(B)$ to that of $\G(A)$, 
so that $\f$ maps the short exact sequence (\ref{exact}) for $\G(B)$ 
to that for $\G(A)$. The above fact also shows that the map $\f$ restricted 
to $\Z^n$ agrees with $P$ for $i, j$ with $B_j^i=0$, in particular, 
if we view the restricted map as a block matrix as before, then it is 
an upper triangular block matrix and the diagonal blocks agree with those of 
$P$. 
Therefore the determinant of the restricted map is equal to $\det P$. This 
proves (1). 

On the other hand, it follows from the definition (\ref{f}) of $\f$ that 
the map induced from $\f$ on $\G(B)/\Z^n=(\Z/2)^n$ is nothing but $P$, 
so it is an isomorphism, proving (2).  
These imply that $\f$ is always injective and an isomorphism if and only if 
$\det P=\pm 1$. 
\end{proof}

\section{Group extension}\label{sect6}

\smallskip
A square $(0,1)$-matrix of size $m$ is in $\GL(m;\Z)$ 
if and only if it is in $\GL(m;\Z/2)$ when $m\le 3$. 
Therefore, if $n_i\le 3$ for all $i$, where $(n_1,\dots,n_q)$ is the type 
of $A$ and $B$, then $\det P=\pm 1$ and $\f$ in Lemma~\ref{odd} 
is an isomorphism. In general $\f$ may not be an isomorphism, but we prove 
the following using the existence of $\g$. 

\begin{lemm} \label{iso}
%$\f(G(A))$ is isomorphic to $G(B)$, so 
If $H^*(M(A);\Z/2)$ is isomorphic to $H^*(M(B);\Z/2)$ as graded rings, then 
$\G(A)$ is isomorphic to $\G(B)$. 
\end{lemm}

We admit the lemma above for the moment and complete the proof of 
Theorem~\ref{main} in the Introduction. 

\begin{proof}[Proof of Theorem~\ref{main}]
Real Bott manifolds are compact riemannian flat manifolds by 
Lemma~\ref{rflat}, hence by a theorem of Bieberbach 
they are diffeomorphic if and only if their fundamental groups 
are isomorphic (see \cite{wolf77}, Theorem 3.3.1 in p.105). 
Therefore, Theorem~\ref{main} follows from Lemma~\ref{iso}. 
\end{proof}

The rest of this section is devoted to the proof of Lemma~\ref{iso}. 
Remember the group extension (\ref{exact}) 
\begin{equation*}\label{ext1}
0\ra \ZZ^n\ra \G(A)\ra (\ZZ_2)^n\ra 1.
\end{equation*}
Conjugation action of $\G(A)$ on $\ZZ^n$ induces a homomorphism 
%a $(\ZZ_2)^n$-action on $\ZZ^n$
\[
\phi_A:(\ZZ_2)^n\ra \mathop{\Aut}(\ZZ^n).
\]
We remark that the $(\ZZ_2)^n$-module $\ZZ^n$ via 
$\phi_A$ decomposes into 
sum of rank one $(\ZZ_2)^n$-modules, which follows from Lemma~\ref{lemm1}. 
There is a $2$-cocycle
\[
f_A:(\ZZ_2)^n\times (\ZZ_2)^n\ra \ZZ^n \] whose cohomology class
$[f_A]\in H^2_{\A}((\ZZ_2)^n;\ZZ^n)$ represents the above group extension, 
that is, $\G(A)$ is the product
$\ZZ^n\times (\ZZ_2)^n$ with group law:
\begin{equation}\label{productA}
(\n,\al)(m,\be)=(\n+\phi_A(\al)(m)+f_A(\al,\be), \al\be).
\end{equation}
Similarly we have $\phi_B$ and $f_B$ for the group $\G(B)$. 

Lemma~\ref{odd} shows that there is a commutative diagram:
\begin{equation*}\label{mapA-B}
\begin{CD}
0@>>> \ZZ^n@>>>\G(B)@>>> (\ZZ_2)^n @>>> 1\\
@.     @V\g VV        @V\g VV     @V{\Phi}VV \\
0 @>>> \ZZ^n@>>> \G(A) @>>>(\ZZ_2)^n @>>> 1.
\end{CD}
\end{equation*}
where $\Phi$ is an isomorphism. 
We write $$\g(0,\al)=(\lambda(\al),\Phi(\al)).$$
Then, for $(\n,\al)\in \G(B)$ we have 
\begin{equation} \label{h2}
\begin{split}
\g(\n,\al)&=\g((\n,1)(0,\al))=\g(\n,1)\g(0,\al)\\
&=(\g(\n),1)(\lambda(\al),\Phi(\al))\\
&=(\g(\n)+\lambda(\al),\Phi(\al)).
\end{split}
\end{equation}
Therefore, applying $\g$ to the both sides of the identity 
$\displaystyle (0,\al)(0,\be)=(f_B(\al,\be), \al\be)$, we have 
\begin{equation*}\begin{split}
\g((0,&\al)(0,\be))=
(\lambda(\al),\Phi(\al))(\lambda(\be),\Phi(\be))\\
&=(\lambda(\al)+\phi_A(\Phi(\al))(\lambda(\be))
+f_A(\Phi(\al),\Phi(\be)), \Phi(\al\be)),
\end{split}
\end{equation*}
while we have 
\[
\g(f_B(\al,\be), \al\be)=
(\g(f_B(\al,\be))+\lambda(\al\be),\Phi(\al\be))
\]
by (\ref{h2}). 
It follows that
\begin{equation}\label{coin}
\g(f_B(\al,\be))=\lambda(\al)+\phi_A(\Phi(\al))(\lambda(\be))-\lambda(\al\be)
+\Phi^*f_A(\al,\be).
\end{equation}

Similarly, applying $\g$ to the both sides of the identity 
$(0,\al)(\n,1)=(\phi_B(\al)(\n),\al)$, we have 
\[
\begin{split}
\g((0,\al)(\n,1))&=(\lambda(\al),\Phi(\al))(\g(\n),1)\\
&=(\lambda(\al)+\phi_A(\Phi(\al))(\g(\n)),\Phi(\al)),
\end{split}
\]
while 
\[
\g(\phi_B(\al)(\n),\al)=(\g(\phi_B(\al)(\n))+\lambda(\al),\Phi(\al)). 
\] 
It follows that
\begin{equation*}\label{commutativity}
\g(\phi_B(\al)(\n))=\phi_A(\Phi(\al))(\g(\n)).
\end{equation*}

We regard elements in $\ZZ^n$ as column vectors and 
represent the homomorphism $\g:\ZZ^n\ra \ZZ^n$ by an integral matrix $\Q$. 
%with respect to a row basis. 
Then the identity above is equivalent to 
\begin{equation*}\label{0st}
\Q\cdot\phi_B(\al)=\phi_A(\Phi(\al))\cdot\Q.
\end{equation*}
We note that $\tQ=(\det Q)Q^{-1}$ is an integral matrix, where $\det\Q$ 
that is the order of $\ZZ^n/\g(\ZZ^n)$ is odd by Lemma~\ref{odd}.  
It follows from the identity above that 
\begin{equation}\label{commutativityinverse}
\phi_B(\al)\cdot\tQ=\tQ\cdot\phi_A(\Phi(\al)).
\end{equation}

Applying $\tQ$ to the both sides of \eqref{coin}, 
we have 
\begin{equation*} \label{dQ}
\begin{split}
&\tQ\Q f_B(\al,\be)\\
=&\tQ\lambda(\al)+\tQ\phi_A(\Phi(\al))(\lambda(\be))
-\tQ\lambda(\al\be)+\tQ\Phi^*f_A(\al,\be)\\
=&\tQ\lambda(\al)+\phi_B(\al)(\tQ\lambda(\be))
-\tQ\lambda(\al\be)+\tQ\Phi^*f_A(\al,\be)\\
=&\delta_B(\tQ\lambda)(\al,\be)+\tQ\Phi^*f_A(\al,\be)
\end{split}\end{equation*}
where we used \eqref{commutativityinverse} at the second identity and 
the definition of coboundary $\delta_B$ at the last identity. 
Since $\tQ Q$ is $\det Q$ times the identity matrix, the identity above 
implies that 
\begin{equation}\label{3rd}
\begin{split}
[\det \Q\cdot f_B]=[\tQ\Phi^*f_A]\in H_{\B}^2((\ZZ_2)^n,\tQ\ZZ^n).
\end{split}\end{equation}
Here $\tQ\ZZ^n$ is viewed as a $(\ZZ_2)^n$-module via $\phi_B$. 
It decomposes into the direct sum of rank one $(\ZZ_2)$-modules because 
we have \eqref{commutativityinverse} and the $(\ZZ_2)^n$-module 
$\ZZ^n$ via $\phi_A$ decomposes into the direct sum of rank one 
$(\ZZ_2)$-modules. Therefore 
\begin{equation} \label{split}
H_{\B}^2((\ZZ_2)^n,\tQ\ZZ^n)\cong\bigoplus_{i=1}^nH_{\phi_i}^2((\Z_2)^n;\Z)
\end{equation}
where $\phi_i\colon (\Z_2)^n\to \Aut(\Z)=\{\pm 1\}$ is a homomorphism. 

\medskip
\noindent
{\bf Fact.} $H_{\phi}^2((\Z_2)^n;\Z)$ is an elementary abelian $2$-group 
for any homomorphism $\phi\colon (\Z_2)^n\to \Aut(\Z)$. 

\medskip
\noindent
(This fact is probably known but since we do not know the literature, 
we will give a proof in the Appendix.) It follows from \eqref{split} and 
the fact above that $H_{\B}^2((\ZZ_2)^n,\tQ\ZZ^n)$ is an elementary 
abelian $2$-group. 
Since $\det Q$ is odd as remarked before, the identity \eqref{3rd} implies that
\begin{equation}\label{4th}
[f_B]=[\tQ\Phi^*f_A]\in H_{\B}^2((\ZZ_2)^n,\tQ\ZZ^n).
\end{equation}

The group $\H$ corresponding to the cocycle $\tQ\Phi^*f_A$ 
is the product $\tQ\ZZ^n\times
(\ZZ_2)^n$ with group law:
\begin{equation}\label{productB}
\begin{split}
&(\tQ\n,\al)(\tQ m,\be)\\
=&(\tQ\n+\phi_B(\al)(\tQ m)+ \tQ f_A(\Phi(\al),\Phi(\be)), \al\be)
\end{split}\end{equation}
in which we note that $\tQ\n+\phi_B(\al)(\tQ m)+
\tQ f_A(\Phi(\al),\Phi(\be))\in \tQ\ZZ^n$. In fact, using
\eqref{commutativityinverse}
\begin{equation}\label{productC}
\begin{split}
&\tQ\n+\phi_B(\al)(\tQ m)+ \tQ f_A(\Phi(\al),\Phi(\be))\\
=&\tQ\n+\tQ\phi_A(\Phi(\al))(m)+\tQ f_A(\Phi(\al),\Phi(\be))\\
=&\tQ(\n+\phi_A(\Phi(\al))(m)+f_A(\Phi(\al),\Phi(\be))).
\end{split}\end{equation}
Since $\G(B)$ is isomorphic to $\H$ by \eqref{4th}, 
it suffices to prove that $\H$ is isomorphic to $\G(A)$.

Define a map $\mathcal T:\G(A)\ra \H$ by 
\[
\mathcal T(\n,\al)=(\tQ\n,\Phi^{-1}(\al)).
\]
This is clearly a bijection. 
Using \eqref{productB}, \eqref{productC} and \eqref{productA}, we have 
\begin{equation*}\label{productD}
\begin{split}
&\mathcal T(\n,\al)\mathcal T(m,\be)\\
=&(\tQ\n,\Phi^{-1}(\al))(\tQ m,\Phi^{-1}(\be))\\
=&(\tQ(\n+\phi_A(\al)(m)+f_A(\al,\be)),\Phi^{-1}(\al\be))\\
=&\mathcal T((\n+\phi_A(\al)(m)+f_A(\al,\be),\al\be))\\
=&\mathcal T((\n,\al)(m,\be)).
\end{split}\end{equation*}
Hence $\mathcal T$ is an isomorphism of $\G(A)$ onto $\H$. 
This completes the proof of Lemma~\ref{iso}.

\begin{rema}
The above proof of Lemma~\ref{iso} actually proves that any subgroup of 
$\G(A)$ with odd index is isomorphic to $\G(A)$. 
\end{rema}

\section{Classification of real Bott manifolds of low dimension} \label{sect7}

%We give some explicit diffeomorphism classification results on real Bott 
%manifolds in dimensions up to 4. 

Real Bott manifolds are determined by upper triangular 
square $(0,1)$-matrices with zero diagonal entries and 
the diffeomorphism classification of real Bott manifolds 
reduces to the isomorphism classification of associated cohomology 
rings with $\Z/2$ coefficients by our main Theorem~\ref{main}. 
As observed in Section~\ref{coho}, we may assume that our matrices are 
of the form \eqref{block} which we call a \emph{normal} form. 
Therefore, it suffices to check which matrices of normal form produce 
isomorphic cohomology rings and this can be done by an elementary computation 
when the size $n$ 
of matrices, that is the dimension of real Bott manifolds, is up to $4$. 
We remember that permuting the suffixes of the cohomology 
generators $x_1,\dots,x_n$ in Section~\ref{coho} 
corresponds to conjugating our matrices by a permutation matrix. 
So the 
cohomology rings associated with conjugate matrices by permutation matrices 
are isomorphic. This decreases necessary computations. 
Below are the results. The same results are obtained in \cite{nazr08} 
independently but the method is different from ours. 
%The authors are also informed that Suyoung Choi \cite{choi08-2} 
%obtained the same result. 

\medskip
\noindent
{\bf The case $n=2$.}  Real Bott manifolds of dimension 2 are 
the torus $(S^1)^2$ 
or the Klein bottle 
and the corresponding matrices of normal form 
are respectively the zero matrix of size 2 and 
{\tiny 
$\begin{pmatrix} 0&1\\
0&0\end{pmatrix}$}. 

\medskip
\noindent
{\bf The case $n=3$.} There are four diffeomorphism classes in real Bott 
manifolds of dimension 3 and corresponding matrices of normal form 
are distinguished by their types as seen below.  
The number of an item below with $\star$ shows that the corresponding real 
Bott manifold is orientable (see Lemma~\ref{orien}). 
%The last two real Bott manifolds below are indecomposable. 
\begin{enumerate}
\item[$1^\star$.] Type $(3)$  

The zero matrix of size 3 
and the real Bott manifold is $(S^1)^3$. 
\item[2.\ ] Type $(2,1)$ 
{\tiny 
\[
\begin{pmatrix} 
0 & 0 & 0\\
0 & 0 & 1\\
0 & 0 & 0\end{pmatrix} \quad \begin{pmatrix} 
0 & 0 & 1\\
0 & 0 & 0\\
0 & 0 & 0\end{pmatrix}\quad \begin{pmatrix} 
0 & 0 & 1\\
0 & 0 & 1\\
0 & 0 & 0\end{pmatrix}
\]}
The real Bott manifold is 
$S^1\times \text{(Klein bottle)}$. 
\item[$3^\star$.] Type $(1,2)$
{\tiny 
\[
\begin{pmatrix} 
0 & 1 & 1\\
0 & 0 & 0\\
0 & 0 & 0\end{pmatrix}
\]}
\item[4.\ ] Type $(1,1,1)$
{\tiny
\[
\begin{pmatrix} 
0 & 1 & 0\\
0 & 0 & 1\\
0 & 0 & 0\end{pmatrix} \quad \begin{pmatrix} 
0 & 1 & 1\\
0 & 0 & 1\\
0 & 0 & 0\end{pmatrix}
\]}
\end{enumerate} 

\begin{rema}
Compact riemannian flat manifolds of dimension 3 are classified.  
There are ten diffeomorphism classes and six of them are orientable 
(\cite[p.117 and p.120]{wolf77}).  One can easily check that the real Bott 
manifolds 
in $1^\star$ and $3^\star$ above are respectively of types 
$\mathcal G_1$ and $\mathcal G_2$ in \cite[Theorem 3.5.5]{wolf77} and 
those in 2 and 4 above are respectively 
of types $\mathcal B_1$ and $\mathcal B_3$ in \cite[Theorem 3.5.9]{wolf77}.  
\end{rema}

\noindent
{\bf The case $n=4$.} There are twelve diffeomorphism classes in real Bott 
manifolds of dimension 4 and corresponding matrices of normal form 
are as described below.  We list representatives of conjugacy classes in 
matrices of normal form by permutation matrices. 
The suffix of a matrix below denotes the number of elements in the 
conjugacy class represented by the matrix. 
The number of an item below with $\star$ shows that the corresponding real 
Bott manifold is orientable as before.  
%The last seven real Bott manifolds below are indecomposable. 
\begin{enumerate}
\item[$1^\star$.] Type $(4)$  

The zero matrix of size 4 and the real Bott manifold is $(S^1)^4$. 
\item[2.\ ] Type $(3,1)$
{\tiny 
\[
\begin{pmatrix} 
0 & 0 & 0 & 0\\
0 & 0 & 0 & 0\\
0 & 0 & 0 & 1\\
0 & 0 & 0 & 0 \end{pmatrix}_3 \quad 
\begin{pmatrix} 
0 & 0 & 0 & 0\\
0 & 0 & 0 & 1\\
0 & 0 & 0 & 1\\
0 & 0 & 0 & 0 \end{pmatrix}_3 \quad 
\begin{pmatrix} 
0 & 0 & 0 & 1\\
0 & 0 & 0 & 1\\
0 & 0 & 0 & 1\\
0 & 0 & 0 & 0 \end{pmatrix}_1
\]}
The real Bott manifold is $(S^1)^2\times 
\text{(Klein bottle)}$. 
\item[$3^\star$.] Type $(2,2)$ 
{\tiny 
\[
\begin{pmatrix} 
0 & 0 & 0 & 0\\
0 & 0 & 1 & 1\\
0 & 0 & 0 & 0\\
0 & 0 & 0 & 0 \end{pmatrix}_2 \quad 
\begin{pmatrix} 
0 & 0 & 1 & 1\\
0 & 0 & 1 & 1\\
0 & 0 & 0 & 0\\
0 & 0 & 0 & 0 \end{pmatrix}_1 
\]}
The real Bott manifold is the product of $S^1$ and 
the 3-dimensional real Bott manifold of Type $(1,2)$.
\item[4.\ ] Type $(2,2)$
{\tiny 
\[
\begin{pmatrix} 
0 & 0 & 0 & 1\\
0 & 0 & 1 & 0\\
0 & 0 & 0 & 0\\
0 & 0 & 0 & 0 \end{pmatrix}_2 \quad 
\begin{pmatrix} 
0 & 0 & 1 & 0\\
0 & 0 & 1 & 1\\
0 & 0 & 0 & 0\\
0 & 0 & 0 & 0 \end{pmatrix}_4 
\]}
The real Bott manifold is (Klein bottle)$\times$(Klein bottle). 
\item[5.\ ] Type $(2,1,1)$ 
{\tiny 
\[
\begin{pmatrix} 
0 & 0 & 0 & 0\\
0 & 0 & 1 & 0\\
0 & 0 & 0 & 1\\
0 & 0 & 0 & 0 \end{pmatrix}_2 \quad 
\begin{pmatrix} 
0 & 0 & 0 & 0\\
0 & 0 & 1 & 1\\
0 & 0 & 0 & 1\\
0 & 0 & 0 & 0 \end{pmatrix}_2 \quad 
\begin{pmatrix} 
0 & 0 & 1 & 0\\
0 & 0 & 1 & 0\\
0 & 0 & 0 & 1\\
0 & 0 & 0 & 0 \end{pmatrix}_1 \quad 
\begin{pmatrix} 
0 & 0 & 1 & 1\\
0 & 0 & 1 & 1\\
0 & 0 & 0 & 1\\
0 & 0 & 0 & 0 \end{pmatrix}_1 
\]}
The real Bott manifold is the product of $S^1$ and 
the 3-dimensional real Bott manifold of Type $(1,1,1)$.
\item[6.\ ] Type $(2,1,1)$ 
{\tiny 
\[
\begin{pmatrix} 
0 & 0 & 0 & 1\\
0 & 0 & 1 & 0\\
0 & 0 & 0 & 1\\
0 & 0 & 0 & 0 \end{pmatrix}_2 \quad 
\begin{pmatrix} 
0 & 0 & 1 & 0\\
0 & 0 & 1 & 1\\
0 & 0 & 0 & 1\\
0 & 0 & 0 & 0 \end{pmatrix}_2 \quad 
\begin{pmatrix} 
0 & 0 & 0 & 1\\
0 & 0 & 1 & 1\\
0 & 0 & 0 & 1\\
0 & 0 & 0 & 0 \end{pmatrix}_2 
\]}
\item[7.\ ] Type $(1,3)$
{\tiny 
\[
\begin{pmatrix} 
0 & 1 & 1 & 1\\
0 & 0 & 0 & 0\\
0 & 0 & 0 & 0\\
0 & 0 & 0 & 0 \end{pmatrix}_1
\]}
\item[8.\ ] Type $(1,2,1)$ 
{\tiny 
\[
\begin{pmatrix} 
0 & 1 & 1 & 0\\
0 & 0 & 0 & 0\\
0 & 0 & 0 & 1\\
0 & 0 & 0 & 0 \end{pmatrix}_2\quad 
\begin{pmatrix} 
0 & 1 & 1 & 0\\
0 & 0 & 0 & 1\\
0 & 0 & 0 & 1\\
0 & 0 & 0 & 0 \end{pmatrix}_1\quad 
\begin{pmatrix} 
0 & 1 & 1 & 1\\
0 & 0 & 0 & 0\\
0 & 0 & 0 & 1\\
0 & 0 & 0 & 0 \end{pmatrix}_2\quad 
\begin{pmatrix} 
0 & 1 & 1 & 1\\
0 & 0 & 0 & 1\\
0 & 0 & 0 & 1\\
0 & 0 & 0 & 0 \end{pmatrix}_1
\]}
\item[$9^\star$.] Type $(1,1,2)$
{\tiny 
\[
\begin{pmatrix} 
0 & 1 & 1 & 0\\
0 & 0 & 1 & 1\\
0 & 0 & 0 & 0\\
0 & 0 & 0 & 0 \end{pmatrix}_2
\]}
\item[10.\ ] Type $(1,1,2)$
{\tiny 
\[
\begin{pmatrix} 
0 & 1 & 0 & 0\\
0 & 0 & 1 & 1\\
0 & 0 & 0 & 0\\
0 & 0 & 0 & 0 \end{pmatrix}_1 \quad
\begin{pmatrix} 
0 & 1 & 1 & 1\\
0 & 0 & 1 & 1\\
0 & 0 & 0 & 0\\
0 & 0 & 0 & 0 \end{pmatrix}_1
\]}
\item[11.\ ] Type $(1,1,1,1)$ 
{\tiny 
\[
\begin{pmatrix} 
0 & 1 & 0 & 1\\
0 & 0 & 1 & 0\\
0 & 0 & 0 & 1\\
0 & 0 & 0 & 0 \end{pmatrix}_1 \quad
\begin{pmatrix} 
0 & 1 & 1 & 0\\
0 & 0 & 1 & 1\\
0 & 0 & 0 & 1\\
0 & 0 & 0 & 0 \end{pmatrix}_1 \quad 
\begin{pmatrix} 
0 & 1 & 1 & 1\\
0 & 0 & 1 & 0\\
0 & 0 & 0 & 1\\
0 & 0 & 0 & 0 \end{pmatrix}_1 
\]}
\item[12.\ ] Type $(1,1,1,1)$ 
{\tiny 
\[
\begin{pmatrix} 
0 & 1 & 0 & 0\\
0 & 0 & 1 & 0\\
0 & 0 & 0 & 1\\
0 & 0 & 0 & 0 \end{pmatrix}_1 \quad
\begin{pmatrix} 
0 & 1 & 1 & 0\\
0 & 0 & 1 & 0\\
0 & 0 & 0 & 1\\
0 & 0 & 0 & 0 \end{pmatrix}_1 \quad 
\begin{pmatrix} 
0 & 1 & 0 & 0\\
0 & 0 & 1 & 1\\
0 & 0 & 0 & 1\\
0 & 0 & 0 & 0 \end{pmatrix}_1 \quad 
\begin{pmatrix} 
0 & 1 & 1 & 1\\
0 & 0 & 1 & 1\\
0 & 0 & 0 & 1\\
0 & 0 & 0 & 0 \end{pmatrix}_1 
\]}
\end{enumerate}

\medskip

\section{Riemannian flat real toric manifolds} \label{sect8}

A toric manifold $X$ of complex dimension $n$ supports an action of $(\C^*)^n$ 
and its real part $\XR$ supports an action of $(\R^*)^n$, where 
$\C^*=\C\backslash\{0\}$ and $\R^*=\R\backslash\{0\}$. 
Let $T$ be the maximal compact toral subgroup of $(\C^*)^n$ and $\TR$ be the 
maximal elementary abelian $2$-subgroup of $(\R^*)^n$. 
The orbit space $\XR/\TR$ can naturally be identified with $X/T$. 
When $X$ is projective, the orbit 
space can be identified with a simple $n$-polytope via a moment map.  

The action of $\TR$ on the real Bott manifold $M(A)=\R^n/\G(A)$ is given 
as follows. Let $r_j$ $(j=1,\dots,n)$ be an involution on $\R^n$ defined by 
\[
r_j(x^1,\dots,x^n)=(x^1,\dots,x^{j-1}, -x^j, x^{j+1},\dots, x^n).
\]
As easily checked 
\[
r_js_i=\begin{cases} s_ir_j \quad&\text{if $i\not=j$,}\\
s_i^{-1}r_j \quad&\text{if $i=j$,}
\end{cases}
\]
where $s_i$ is the Euclidean motion on $\R^n$ defined in \eqref{si}, so 
$r_j$ induces an involution $\r_j$ on $M(A)=\R^n/\G(A)$. Obviously $\r_j$'s 
commute with each other so that they generate an elementary 
abelian $2$-group of rank $n$ and this gives the action of $\TR$. 

We remark that the action of $\TR$ on $M(A)=\R^n/\G(A)$ preserves 
the flat riemannian metric on it. 
The group generated by $s_i$'s and $r_j$'s agrees with the group generated 
by $r_j$'s and translations by $\frac{1}{2}e^1,\dots,\frac{1}{2}e^n$ 
where $e^1,\dots,e^n$ are the standard basis of $\R^n$ 
as before.  It follows that the orbit space $M(A)/\TR$ can be identified 
with an $n$-cube 
\[
\{(x^1,\dots,x^n)\in\R^n \mid 0\le x^1\le 1/2,\dots, 0\le x^n\le 1/2\}.
\] 

The purpose of this section is to prove Theorem~\ref{main2} in the 
Introduction, that is

\begin{theo} \label{flat}
A real toric manifold of dimension $n$ which 
admits a flat riemannian metric invariant under the action of $\TR$ 
is a real Bott manifold. 
\end{theo}

We recall some results for the proof of the theorem above. 
Let $X$ be a toric manifold and let $X_i$ $(1\le i\le m)$ be a  
connected complex codimension-one closed submanifold of $X$ 
fixed pointwise under some circle subgroup $T_i$ of the torus $T$.  
We call $X_i$ a \emph{characteristic submanifold} of $X$. 
Then 
\[
K_X:=\{ I\subset \{1,\dots,m\}\mid \cap_{i\in I}X_i\not=\emptyset\}
\]
is the underlying abstract simplicial complex of the fan of $X$. 

Let $\XR$ be the real part of $X$.  
The intersection $X_i\cap \XR$ is a connected real codimension-one 
closed submanifold of $\XR$ fixed pointwise under the order two subgroup 
$T_i\cap\TR$ of $\TR$.  
Conversely any connected real codimension-one 
closed submanifold of $\XR$ fixed pointwise under an order two subgroup 
of $\TR$ is the intersection of $\XR$ with some $X_i$. 
We call those closed submanifolds \emph{characteristic submanifolds} of 
$\XR$ as well.  
This observation says that there is a bijective correspondence 
between characteristic submanifolds of $X$ and those of $\XR$.  
Hence one can also define $K_X$ using the characteristic submanifolds 
of $\XR$. 

We say that a simplicial complex is a \emph{crosscomplex} of dimension $n-1$ 
if it is the boundary complex of a crosspolytope of dimension $n$, where 
a crosspolytope of dimension $n$ is the dual (or polar) of an $n$-cube. 
We recall two facts from \cite{ma-pa08}.  The first lemma below is stated 
in \cite[Corollary 3.5]{ma-pa08} in the complex case but it also holds 
in the real case as stated because of the observation above. 

\begin{lemm}[Corollary 3.5 in \cite{ma-pa08}]  \label{KX}
A real toric manifold $\XR$ is a real Bott manifold if and only if 
the simplicial complex $K_X$ associated with $\XR$ is a crosscomplex. 
\end{lemm}

\begin{lemm}[Lemma 4.7 in \cite{ma-pa08}] \label{cross}
Let $K$ be a connected simplicial complex of dimension $k\ge 2$. 
If the link of each vertex of $K$ is a crosscomplex of dimension $k-1$, 
then $K$ is a crosscomplex. 
\end{lemm}

\begin{proof}[Proof of Theorem~\ref{flat}]
We shall prove the theorem by induction on the dimension $n$. 
The theorem is obvious when $n=1$. 
A closed surface which admits a flat riemannian metric is a torus or 
a Klein bottle and they are real Bott manifolds, so the theorem also holds 
when $n=2$.  

Now suppose the theorem holds for $n-1\ge 2$ and 
let $\XR$ be a real toric manifold of dimension $n$ which 
satisfies the assumption in the theorem.  
Let $\XR_1,\dots,\XR_m$ be the characteristic submanifolds of $\XR$. 
A vertex of the simplicial complex $K_X$ associated with $\XR$ 
corresponds to some $\XR_i$ 
and the link of the vertex is the simplicial complex associated with $\XR_i$. 
Since $\XR$ admits 
a riemannian flat metric invariant under the action of $\TR$, 
each $\XR_i$ is again a riemannian flat manifold  
because it is fixed pointwise under a subgroup of $\TR$. 
Therefore $\XR_i$ is a real Bott manifold 
by the inductive assumption and hence the link of the vertex of $K_X$ 
is crosscomplex by Lemma~\ref{KX}.  
Since $\dim K_X=n-1\ge 2$, $K_X$ is a crosscomplex 
by Lemma~\ref{cross} and hence $X$ is a real Bott manifold 
by Lemma~\ref{KX}. This completes the 
induction step and the proof of the theorem. 
\end{proof}

\section{Small cover} \label{sect9}

Let $\TR$ be an elementary abelian $2$-group of rank $n$ as before. 
A closed smooth manifold $M$ of dimension $n$ with a smooth action of $\TR$ 
is called locally standard if each point of $M$ has an invariant open 
neighborhood equivariantly diffeomorphic to an invariant open subset of 
a faithful real $\TR$-module of dimension $n$. 
The orbit space of a locally standard $\TR$-manifold $M$ 
is a manifold with corners 
because the orbit space of a faithful $\TR$-module of dimension $n$ 
is homeomorphic to the product of $n$ half lines.  
A convex polytope of dimension $n$ is called simple if there are exactly 
$n$ edges meeting at each vertex, and a simplex convex polytope is a typical 
example of a manifold with corners.  
If $M$ is locally standard and the orbit space is identified with 
a simple convex 
polytope $P$, then $M$ is called a \emph{small cover} over $P$ 
(\cite{da-ja91}).  

A real toric manifold $\XR$ with the natural 
$\TR$-action is locally standard and its orbit space is often a simple 
convex polytope. In fact, this is the case when $X$ is projective, so 
a real toric manifold $\XR$ is a small cover when $X$ is projective. 
However there are many small covers which do not arise this way. 
For example, every closed surface becomes a small cover 
but only the torus $S^1\times S^1$ is a real toric 
manifold among orientable closed surfaces (\cite{payn04}).  
We may think of small covers as a topological counterpart to 
real toric manifolds and may ask the same question as in the 
Introduction for small covers. We remark that equivariant homeomorphism types 
of small covers can be distinguished by their equivariant cohomology algebras with $\Z/2$ 
coefficients (\cite{masu07}). 

When $\XR$ is a real Bott manifold, the orbit space is an $n$-cube 
as observed in Section~\ref{sect8}; so a real Bott manifold of dimension $n$ 
becomes a small cover over an $n$-cube and 
the converse is known to be true up to homeomorphism. 

\begin{theo}[\cite{ma-pa08}, \cite{ch-ma-su08-1}]
A small cover over an $n$-cube is homeomorphic to a real Bott manifold of 
dimension $n$. 
\end{theo}

The number $Q_n$ of equivariant homeomorphism classes in small covers over 
an $n$-cube is computed in \cite{choi08} for any $n$, e.g. 
$$Q_1=1,\ Q_2=6,\ Q_3=259,\ Q_4=87360,\ Q_5=236240088, \dots.$$
However, the number $H_n$ of (non-equivariant) 
homeomorphism classes in small covers over 
an $n$-cube is unknown although 
\[
H_1=1,\ H_2=2,\ H_3=4,\ H_4=12
\]
as described in Section~\ref{sect7}.  
%One can also see that
%$H_n>2^{(n-2)(n-3)/2}$ for any $n$ 
%(\cite{masu08}). 

%and $H_n\ge 2^{n-1}$ for any $n$ because the number of types of cohomology 
%rings of small covers over an $n$-cube is $2^{n-1}$ and there exists at 
%least one small cover for each type. 

As is well-known, regular simple polytopes of dimension $n\ge 3$ 
are an $n$-cube and an $n$-simplex in each dimension $n$, 
the dodecahedron in dimension 3 
and the 120-cell in dimension 4.  The homeomorphism type of a small cover 
over an $n$-simplex is unique, that is the real projective space of dimension 
$n$. Small covers over the dodecahedron and the 120-cell admit hyperbolic 
metrics and are studied in \cite{ga-sc02} from this point of view. 
In particular, it is proved 
in the paper that there are exactly $25$ small covers over the 
dodecahedron up to isometry (equivalently up to homeomorphism by Mostow 
rigidity).

\section*{Appendix}

In this appendix we give a proof of the Fact used in Section~\ref{sect6}. 
In fact we will prove a more precise statement. 
It is well-known that $H_{\phi}^2((\Z_2)^n;\Z)$ is isomorphic to 
$(\Z/2)^{n}$ when $\phi$ is trivial. We prove 

\begin{theoap}
\emph{If $\phi$ is non-trivial, 
$H_{\phi}^2((\Z_2)^n;\Z)$ is isomorphic to $(\Z/2)^{n-1}$.}
\end{theoap}

We recall the following Hochschild-Serre spectral sequence, 
see \cite[p.355]{Mac} or \cite{HS}. 

\begin{propap}\label{HoSe}
\emph{Let $1\ra \Gamma\ra \Pi\ra
\Pi/\Gamma\ra 1$ be a group extension and let $A$ be a $\Pi$-module
through a homomorphism $\phi: \Pi\ra \mathop{\Aut}(A)$.
Suppose $m\geq 1$ and $H^q_{\phi}(\Gamma,A)=0$ for $1<q<m$. For
$0<q<m$, there is the exact sequence
\begin{equation*}\label{Hoch-Se}
\tag{A.1}
\begin{split}
&0\ra H^1_{\phi}(\Pi/\Gamma,A^{\Gamma})\ra H^1_{\phi}(\Pi,A)\ra
H^0_{\phi}(\Pi/\Gamma,H^1_{\phi}(\Gamma,A))\ra\\
&\cdots\ra H^q_{\phi}(\Pi/\Gamma,A^{\Gamma})\ra H^q_{\phi}(\Pi,A)\ra
H^{q-1}_{\phi}(\Pi/\Gamma,H^1_{\phi}(\Gamma,A))\\
&\ \ \ \ \ \ra H^{q+1}_{\phi}(\Pi/\Gamma,A^{\Gamma})\ra
H^{q+1}_\phi(\Pi,A)\ra\cdots \end{split}
\end{equation*}}
%\begin{equation*}\label{Hoch-Se}
%\tag{A.1}
%\begin{split}
%&\cdots\ra H^q_{\phi}(\Pi/\Gamma,A^{\Gamma})\ra H^q_{\phi}(\Pi,A)\ra
%H^{q-1}_{\phi}(\Pi/\Gamma,H^1_{\phi}(\Gamma,A))\\
%&\ \ \ \ \ \ra H^{q+1}_{\phi}(\Pi/\Gamma,A^{\Gamma})\ra
%H^{q+1}_\phi(\Pi,A)\ra\cdots \end{split}
%\end{equation*}}
\end{propap}

We take $\Pi=(\ZZ_2)^n$ $(n\geq 2)$ and $A=\ZZ$ as a $\Pi$-module through
$\phi: \Pi\ra \mathop{\Aut}(\ZZ)=\{\pm 1\}$. Choose an order two subgroup 
$\Gamma\subset (\ZZ_2)^n$ such that $\phi(\Gamma)=\{\pm 1\}$. 
Clearly 
\begin{equation*}\label{action4}
 \Pi=\Gamma\times \Ker\phi.
\end{equation*}
It is known and easy to check that $H^{2}_\phi(\Gamma,A)=0$, so 
the assumption in the proposition above is 
satisfied for $m=3$. As $A^{\Gamma}=0$ by our condition, 
$H^r_{\phi}(\Pi/\Gamma,A^{\Gamma})=0$ for any $r\geq 0$.
Then the exact sequence \eqref{Hoch-Se} becomes
\begin{equation}\label{isoH}
\tag{A.2} 0\ra H^2_{\phi}(\Pi,A)\ra
H^{1}_{\phi}(\Pi/\Gamma,H^1_{\phi}(\Gamma,A))\ra 0. 
\end{equation}

On the other hand, it is also known and easy to check that 
$H^{1}_\phi(\Gamma,A)\cong \ZZ/2$, so the action of 
$\Pi/\Gamma$ on $H^{1}_\phi(\Gamma,A)$ must be trivial. 
It follows that 
\begin{equation*} \begin{split}
H^{1}_{\phi}(\Pi/\Gamma,H^1_{\phi}(\Gamma,A))&\cong 
H^{1}(\Pi/\Gamma,\ZZ/2)\\
&\cong  H^{1}((\ZZ_2)^{n-1},\ZZ/2)\\
&\cong (\ZZ/2)^{n-1}. \end{split}
\end{equation*} 
This together with \eqref{isoH} implies the theorem. 

\bigskip

\noindent
{\bf Acknowledgment}. This work was motivated by a talk of Suyoung Choi 
on \cite{choi08} 
given at Fudan University in January 2008 and very much stimulated by 
the discussion with Taras Panov which the second author had during his 
stay at Fudan University.  He would like to thank them and also Zhi L\"u 
for inviting him to Fudan University and organizing fruitful seminars. 
He also would like to thank Dong Youp Suh for comments on an 
earlier version of this paper 
and Yasuzo Nishimura for showing him a note concerning matrices 
discussed in Section~\ref{sect7}.

\end{document}